\documentclass{article}%
\usepackage{graphicx}
\usepackage{amsmath}
\usepackage{amsfonts}
\usepackage{amssymb}%
\begin{document}

\title{An Elementary Proof of Hopf's Curvatura Integra Theorem}
\author{Pavel Greenfield\\Drexel University}
\maketitle

\begin{abstract}
We provide an elementary proof of the key aspect of Hopf's curvatura integra
theorem. Namely, we show that the total curvature, i.e. the surface integral
of the Gauss-Kronecker curvature $B$, is independent of the shape of the
hypersurface. We accomplish this task by using the Calculus of Moving Surfaces
to show that the rate of change of the total curvature under smooth changes in
shape vanishes.

\end{abstract}

\section{Introduction}

The Gauss-Bonnet theorem was born out of an elementary observation:\ as a
circle or a sphere gets smaller, its curvature grows larger proportionally.
This gives impetus to the idea that for an arbitrary closed surface, there may
exist a measure related to local curvature whose integral is independent of shape.

This is certainly our experience with curves in two dimensions. When walking
around Central Park in New York, our path may turn this way and that, but when
we complete the loop once, we will have turned a total of $360^{\circ}$ --
unless we also walked around a fountain once in which case we will have turned
a total of $720^{\circ}$. In other words, the total turning equals
$360^{\circ}$ times the number of loops that we have made. Of course, for this
idea to truly work, we must supply our turning measure a sign convention:
counterclockwise loops contribute positively, clockwise loops negatively.

In Vector Calculus, this common experience for any smooth curve $L$ is
expressed precisely by the line integral identity%
\begin{equation}
\int_{L}\kappa dL=2\pi n,
\end{equation}
where $\kappa$ is the signed curvature of the curve and $n$ is the
\textit{turning number}, i.e. the signed number of loops that the curve makes.
Crucially, the integral is independent of the shape of the curve.

For a surface $S$ in three dimensions, the analogous role is played by the
Gauss-Bonnet theorem%
\begin{equation}
\int_{S}KdS=4\pi d, \label{GB}%
\end{equation}
where $K$ is the Gaussian curvature and $d$ is the \textit{degree} of the
surface, a measure analogous to the turning number but more elaborate.

Today, the Gauss-Bonnet theorem is most often discussed in it
\textit{intrinsic} form, where the Gaussian curvature $K$ is derived from the
metric tensor or the Christoffel symbol. In this interpretation, it is a
special case of the Chern-Gauss-Bonnet theorem that generalizes it to higher
(even) dimensions and has important applications in Topology. Our focus will
be on its \textit{extrinsic} interpretation where the manifold is treated as
an $n$-dimensional hypersurface in an $\left(  n+1\right)  $-dimensional
space, for which there is a well-defined curvature tensor $B_{\beta}^{\alpha}$
and $K$ is defined as its determinant. We will refer to this quantity as the
\textit{Gauss-Kronecker curvature} and denote it by the symbol $B$. (For a
planar curve, $B$ reduces to the curvature $\kappa$.) According to the Gauss
equations of the surface, the Gaussian and the Gauss-Kronecker curvatures
coincide for a two-dimensional surface embedded in the three-dimensional
space, and that is why they are often used interchangeably.

In his $1925$ paper \cite{Hopf1925CurvaturaIntegra}, Heinz Hopf proved his
celebrated \textit{curvatura integra theorem}, i.e. the total curvature
theorem, which states that for a closed hypersurface in any dimension, the
integral of $B$ is independent of shape and depends only on the degree of the
surface $d$. Specifically,%
\begin{equation}
\int_{S}BdS=dS_{n},
\end{equation}
where $S_{n}=2\pi^{\left(  n+1\right)  /2}\Gamma^{-1}\left(  \frac{n+1}%
{2}\right)  $ is the "surface area" of the unit $n$-dimensional sphere. 

The object of this paper is to offer an elementary proof of the fact that for
any smooth hypersurface embedded in the Euclidean space, the integral%
\begin{equation}
\int_{S}BdS
\end{equation}
is independent of the shape in the narrow sense that for any two such shapes
related by a smooth transformation, the value of the integral is identical.

Our proof will be based on the Calculus of Moving Surfaces \cite{GrinfeldTC},
which is an extension of Tensor Calculus \cite{McConnellTensors} to deforming
manifolds. Our proof is elementary in the sense that it is based strictly on
the application of analytical rules and does not require a deep understanding
of Differential Geometry or unique logical constructs. Thus, the Calculus of
Moving Surfaces exhibits precisely the kind of robustness, where results or
obtained by a clear application of analytical tools, that Euler described as
having the \textit{highest degree of perfection}.

\section{The tensor description of a hypersurface}

\subsection{Objects on a stationary hypersurface}

Every hypersurface $S$ is endowed with the curvature tensor $B_{\beta}%
^{\alpha}$, \cite{McConnellTensors}, \cite[sec. $12.2$]{GrinfeldTC}, which is
also known as the \textit{second fundamental form}, the \textit{shape
operator}, or the \textit{Weingarten map}. The covariant version
$B_{\alpha\beta}$ is symmetric \cite[eq. $12.36$]{GrinfeldTC}
\begin{equation}
B_{\alpha\beta}=B_{\beta\alpha},
\end{equation}
and therefore $B_{\beta}^{\alpha}$ has a full set of eigenvalues known as the
\textit{principal curvatures} and a full set of orthogonal eigenvectors known
as the \textit{principal directions} or \textit{directions of principal
curvature}. Since the determinant of the matrix equals the product of its
eigenvalues, the Gauss-Kronecker curvature $B$ may also be defined as the
product of the principal curvatures.

The trace $B_{\alpha}^{\alpha}$ of $B_{\beta}^{\alpha}$ is known as the
\textit{mean curvature} and will play the central role in the evolution of the
surface integral formula (\ref{dS}). The determinant of $B_{\beta}^{\alpha}$
is, as stated above, the \textit{Gauss-Kronecker curvature}. For an
$n$-dimensional surface, $B$ is given by the explicit equation \cite[eq.
$9.18$]{GrinfeldTC}%
\begin{equation}
B=\frac{1}{n!}\delta_{\beta_{1}\cdots\beta_{n}}^{\alpha_{1}\cdots\alpha_{n}%
}B_{\alpha_{1}}^{\beta_{1}}\cdots B_{\alpha_{n}}^{\beta_{n}},
\end{equation}
where $\delta_{\beta_{1}\cdots\beta_{n}}^{\alpha_{1}\cdots\alpha_{n}}$ is the
Kronecker $\delta$-symbol \cite[sec. $9.4$]{GrinfeldTC}. For future reference,
the $\delta$-symbol and the determinant of any second-order system in general
and $B_{\beta}^{\alpha}$ in particular satisfy the equation \cite[eq. $9.19$%
]{GrinfeldTC}:%
\begin{equation}
\delta_{\beta_{1}\cdots\beta_{n}}^{\alpha_{1}\cdots\alpha_{n}}B_{\alpha_{1}%
}^{\gamma_{1}}\cdots B_{\alpha_{n}}^{\gamma_{n}}=B\delta_{\beta_{1}\cdots
\beta_{n}}^{\gamma_{1}\cdots\gamma_{n}}\label{db = Bd}%
\end{equation}
Finally, the curvature tensor $B_{\beta}^{\alpha}$ satisfies the
\textit{Codazzi equations} \cite[eq. $12.73$]{GrinfeldTC}%
\begin{equation}
\nabla_{\alpha}B_{\beta}^{\gamma}=\nabla_{\beta}B_{\alpha}^{\gamma
}.\label{Codazzi}%
\end{equation}
In other words, the tensor $\nabla_{\alpha}B_{\beta}^{\gamma}$ is symmetric in
its lower indices. Consequently, if it is contracted on both $\alpha$ and
$\beta$ with an object skew-symmetric in the corresponding indices, such as
any two indices in the $\delta$-symbol, the result will vanish. This argument
will be used in the analysis below.

\subsection{The essential elements of the Calculus of Moving Surfaces}

Consider a smooth evolution $S\left(  t\right)  $ of the hypersurface $S$. The
evolution of $S$ is described by the normal velocity field $C$. The key
differential operator on the moving surface is the invariant time derivative
$\dot{\nabla}$. When applied to the curvature tensor, it gives \cite[eq.
$16.68$]{GrinfeldTC}%
\begin{equation}
\dot{\nabla}B_{\beta}^{\alpha}=\nabla^{\alpha}\nabla_{\beta}C+CB_{\gamma
}^{\alpha}B_{\beta}^{\gamma}\label{.B}%
\end{equation}

Our goal is to show that
\begin{equation}
\frac{d}{dt}\int_{S\left(  t\right)  }BdS=0.\label{dIB = 0}%
\end{equation}
Going forward, we will omit the parameter $t$ and write $S$ instead of
$S\left(  t\right)  $ in order to keep the notation uncluttered.

The key formula that governs the evolution of surface integrals is \cite[eq.
$15.56$]{GrinfeldTC}%
\begin{equation}
\frac{d}{dt}\int_{S}FdS=\int_{S}\dot{\nabla}FdS-\int_{S}CFB_{\alpha}^{\alpha
}dS.\label{dS}%
\end{equation}

\subsection{A proof that $\int_{S}\mathbf{N}dS=\mathbf{0}$}

As a brief aside that will demonstrate our overall approach, we will prove
that%
\begin{equation}
\int_{S}\mathbf{N}dS=\mathbf{0,}%
\end{equation}
where $S$ is a smooth two-dimensional closed surface in the three-dimensional
Euclidean space and $\mathbf{N}$ is  the unit normal. A direct proof of this
identity can be found in \cite{Grinfeld2023VectorTensor}. Here, we will prove
it instead by showing that for a smooth evolution $S\left(  t\right)  $ of the
surface $S$, the integral of the normal remains constant, i.e.%
\begin{equation}
\frac{d}{dt}\int_{S}\mathbf{N}dS=\mathbf{0.}\label{dIN = 0}%
\end{equation}

By equation (\ref{dS}), we have%
\begin{equation}
\frac{d}{dt}\int_{S}\mathbf{N}dS=\int_{S}\dot{\nabla}\mathbf{N}dS-\int
_{S}C\mathbf{N}B_{\alpha}^{\alpha}dS
\end{equation}
The derivative $\dot{\nabla}\mathbf{N}$ of the unit normal is given by the
Thomas formula \cite[eq. $60$]{Thomas1957ExtendedCompatibility}, \cite[eq.
$15.48$]{GrinfeldTC}
\begin{equation}
\dot{\nabla}\mathbf{N}=-\mathbf{S}_{\alpha}\nabla^{\alpha}C,
\end{equation}
where $\mathbf{S}_{\alpha}$ is the covariant basis \cite[eq. $10.05$%
]{GrinfeldTC} and $\nabla^{\alpha}$ is the contravariant surface derivative.
Thus,%
\begin{equation}
\frac{d}{dt}\int_{S}\mathbf{N}dS=-\int_{S}\mathbf{S}_{\alpha}\nabla^{\alpha
}CdS-\int_{S}C\mathbf{N}B_{\alpha}^{\alpha}dS.
\end{equation}
By the product rule in the form $fg^{\prime}=\left(  fg\right)  ^{\prime
}-f^{\prime}g$, we observe that
\begin{equation}
\mathbf{S}_{\alpha}\nabla^{\alpha}C=\nabla^{\alpha}\left(  \mathbf{S}_{\alpha
}C\right)  -\nabla^{\alpha}\mathbf{S}_{\alpha}\ C
\end{equation}
and since \cite[eq. $11.16$]{GrinfeldTC}%
\begin{equation}
\nabla^{\alpha}\mathbf{S}_{\alpha}=\mathbf{N}B_{\alpha}^{\alpha},
\end{equation}
we have%
\begin{equation}
\frac{d}{dt}\int_{S}\mathbf{N}dS=-\int_{S}\nabla^{\alpha}\left(
\mathbf{S}_{\alpha}C\right)  dS+\int_{S}C\mathbf{N}B_{\alpha}^{\alpha}%
dS-\int_{S}C\mathbf{N}B_{\alpha}^{\alpha}dS.
\end{equation}
Now, the first integral vanishes by the divergence theorem and the fact that
$S$ does not have a boundary, while the second and the third integrals cancel
each other.

Thus, the proof of equation (\ref{dIN = 0}) is complete, and we have therefore
demonstrated that the integral of the unit normal is independent of shape as
long as the surface evolves smoothly. To prove that the value of the integral
equals $\mathbf{0}$, we must simply imagine that the surface shrinks to a
point and observe that the integral of the unit normal must tend to zero along
with it. This completes the proof.

We now apply the same approach to our central task: showing that the integral
of $B$ is independent of shape.

\section{A proof that $\int_{S}BdS$ is independent of shape}

The time evolution of the integral of $B$ is governed by the equation%
\begin{equation}
\frac{d}{dt}\int_{S}BdS=\int_{S}\dot{\nabla}BdS-\int_{S}CBB_{\alpha}^{\alpha
}dS.\label{dIB}%
\end{equation}
First focus on the derivative $\dot{\nabla}B$. By the product rule,
$\dot{\nabla}B$ can be expressed as a sum of $n$ similar terms
\begin{equation}
\dot{\nabla}B=\frac{1}{n!}\delta_{\beta_{1}\cdots\beta_{n}}^{\alpha_{1}%
\cdots\alpha_{n}}\sum_{i}B_{\alpha_{1}}^{\beta_{1}}%
\mkern-4mu\cdot\mkern-3mu\cdot\mkern-3mu\cdot\mkern-2mu%
\dot{\nabla}B_{\alpha_{i}}^{\beta_{i}}%
\mkern-4mu\cdot\mkern-3mu\cdot\mkern-3mu\cdot\mkern-2mu%
B_{\alpha_{n}}^{\beta_{n}},
\end{equation}
where we used the metrinilic property \cite[eq. $12.73$]{GrinfeldTC} to leave
the $\delta$-symbol unaffected. Substituting equation (\ref{.B}) into the
identity above yields%
\begin{equation}
\dot{\nabla}B=\frac{1}{n!}\delta_{\beta_{1}\cdots\beta_{n}}^{\alpha_{1}%
\cdots\alpha_{n}}\sum_{i}\left(  B_{\alpha_{1}}^{\beta_{1}}%
\mkern-4mu\cdot\mkern-3mu\cdot\mkern-3mu\cdot\mkern-2mu%
\nabla_{\alpha_{i}}\nabla^{\beta_{i}}C%
\mkern-4mu\cdot\mkern-3mu\cdot\mkern-3mu\cdot\mkern-2mu%
B_{\alpha_{n}}^{\beta_{n}}+CB_{\alpha_{1}}^{\beta_{1}}%
\mkern-4mu\cdot\mkern-3mu\cdot\mkern-3mu\cdot\mkern-2mu%
B_{\alpha_{i}}^{\gamma}B_{\gamma}^{\beta_{i}}%
\mkern-4mu\cdot\mkern-3mu\cdot\mkern-3mu\cdot\mkern-2mu%
B_{\alpha_{n}}^{\beta_{n}}\right)  .\label{.B 1}%
\end{equation}
First, consider the second term in each parenthesized sum with the factor
$B_{\gamma}^{\beta_{i}}$ omitted, i.e.%
\begin{equation}
\delta_{\beta_{1}\cdots\beta_{n}}^{\alpha_{1}\cdots\alpha_{n}}B_{\alpha_{1}%
}^{\beta_{1}}%
\mkern-4mu\cdot\mkern-3mu\cdot\mkern-3mu\cdot\mkern-2mu%
B_{\alpha_{i}}^{\gamma}%
\mkern-4mu\cdot\mkern-3mu\cdot\mkern-3mu\cdot\mkern-2mu%
B_{\alpha_{n}}^{\beta_{n}}.
\end{equation}
Observe that each of the indices $\alpha_{1}\cdots\alpha_{n}$ is represented
among the lower indices of the $n$ factors containing the curvature tensor.
Therefore, according to equation (\ref{db = Bd}), we have
\begin{equation}
\delta_{\beta_{1}\cdots\beta_{n}}^{\alpha_{1}\cdots\alpha_{n}}B_{\alpha_{1}%
}^{\beta_{1}}%
\mkern-4mu\cdot\mkern-3mu\cdot\mkern-3mu\cdot\mkern-2mu%
B_{\alpha_{i}}^{\gamma}%
\mkern-4mu\cdot\mkern-3mu\cdot\mkern-3mu\cdot\mkern-2mu%
B_{\alpha_{n}}^{\beta_{n}}=B\delta_{\beta_{1}\cdots\beta_{i}\cdots\beta_{n}%
}^{\beta_{1}\cdots\gamma\cdots\beta_{n}}.
\end{equation}
Since \cite[eq. $9.23$--$9.25$]{GrinfeldTC}
\begin{equation}
\delta_{\beta_{1}\cdots\beta_{i}\cdots\beta_{n}}^{\beta_{1}\cdots\gamma
\cdots\beta_{n}}=\left(  n-1\right)  !
\end{equation}
we find that%
\begin{equation}
\delta_{\beta_{1}\cdots\beta_{n}}^{\alpha_{1}\cdots\alpha_{n}}B_{\alpha_{1}%
}^{\beta_{1}}%
\mkern-4mu\cdot\mkern-3mu\cdot\mkern-3mu\cdot\mkern-2mu%
B_{\alpha_{i}}^{\gamma}%
\mkern-4mu\cdot\mkern-3mu\cdot\mkern-3mu\cdot\mkern-2mu%
B_{\alpha_{n}}^{\beta_{n}}=\left(  n-1\right)  !B\delta_{\beta_{i}}^{\gamma}.
\end{equation}
Therefore, the second term in each parenthesized sum in equation (\ref{.B 1}),
with the factor $B_{\gamma}^{\beta_{i}}$ now included, reduces to%
\begin{equation}
\frac{1}{n}CB\delta_{\beta_{i}}^{\gamma}B_{\gamma}^{\beta_{i}}=\frac{1}%
{n}CBB_{\alpha}^{\alpha},
\end{equation}
which is a value that is identical for all $n$ terms in the sum. Thus, these
terms combine to produce $CBB_{\alpha}^{\alpha}$, which, of course, cancels
the term $-\int_{S}CBB_{\alpha}^{\alpha}dS$ in (\ref{dIB}). Thus, we are left
with%
\begin{equation}
\frac{d}{dt}\int_{S}BdS=\frac{1}{n!}\delta_{\beta_{1}\cdots\beta_{n}}%
^{\alpha_{1}\cdots\alpha_{n}}\sum_{i}\int_{S}\left(  B_{\alpha_{1}}^{\beta
_{1}}%
\mkern-4mu\cdot\mkern-3mu\cdot\mkern-3mu\cdot\mkern-2mu%
\nabla_{\alpha_{i}}\nabla^{\beta_{i}}C%
\mkern-4mu\cdot\mkern-3mu\cdot\mkern-3mu\cdot\mkern-2mu%
B_{\alpha_{n}}^{\beta_{n}}\right)  dS.
\end{equation}

Apply the product rule pattern $f^{\prime}g=\left(  fg\right)  ^{\prime
}-fg^{\prime}$ to the covariant derivative $\nabla_{\alpha_{i}}$, i.e.%
\begin{align}
\delta_{\beta_{1}\cdots\beta_{n}}^{\alpha_{1}\cdots\alpha_{n}}B_{\alpha_{1}%
}^{\beta_{1}}%
\mkern-4mu\cdot\mkern-3mu\cdot\mkern-3mu\cdot\mkern-2mu%
\nabla_{\alpha_{i}}\nabla^{\beta_{i}}C%
\mkern-4mu\cdot\mkern-3mu\cdot\mkern-3mu\cdot\mkern-2mu%
B_{\alpha_{n}}^{\beta_{n}}  & =\nabla_{\alpha_{i}}\left(  \delta_{\beta
_{1}\cdots\beta_{n}}^{\alpha_{1}\cdots\alpha_{n}}B_{\alpha_{1}}^{\beta_{1}}%
\mkern-4mu\cdot\mkern-3mu\cdot\mkern-3mu\cdot\mkern-2mu%
\nabla^{\beta_{i}}C%
\mkern-4mu\cdot\mkern-3mu\cdot\mkern-3mu\cdot\mkern-2mu%
B_{\alpha_{n}}^{\beta_{n}}\right)  \label{dBCB}\\
& \ \ \ \ \ \ -\delta_{\beta_{1}\cdots\beta_{n}}^{\alpha_{1}\cdots\alpha_{n}%
}\nabla^{\beta_{i}}C\nabla_{\alpha_{i}}\left(  B_{\alpha_{1}}^{\beta_{1}}%
\mkern-4mu\cdot\mkern-3mu\cdot\mkern-3mu\cdot\mkern-2mu%
.%
\mkern-4mu\cdot\mkern-3mu\cdot\mkern-3mu\cdot\mkern-2mu%
B_{\alpha_{n}}^{\beta_{n}}\right)  ,\nonumber
\end{align}
where we have used the metrinilic property of the covariant derivative
\cite[eq. $11.4$]{GrinfeldTC} to leave the $\delta$-symbol outside of the
covariant derivative. When integrated, the first term on the right of the
above equation vanishes by the divergence theorem, combined with the fact that
the surface $S$ does not have a boundary.

In order to proceed with the remaining term%
\begin{equation}
\delta_{\beta_{1}\cdots\beta_{n}}^{\alpha_{1}\cdots\alpha_{n}}\nabla
^{\beta_{i}}C\nabla_{\alpha_{i}}\left(  B_{\alpha_{1}}^{\beta_{1}}%
\mkern-4mu\cdot\mkern-3mu\cdot\mkern-3mu\cdot\mkern-2mu%
.%
\mkern-4mu\cdot\mkern-3mu\cdot\mkern-3mu\cdot\mkern-2mu%
B_{\alpha_{n}}^{\beta_{n}}\right)  ,
\end{equation}
apply the product rule to the covariant derivative $\nabla_{\alpha_{i}}$:%
\begin{align}
& \delta_{\beta_{1}\cdots\beta_{n}}^{\alpha_{1}\cdots\alpha_{n}}\nabla
^{\beta_{i}}C\nabla_{\alpha_{i}}\left(  B_{\alpha_{1}}^{\beta_{1}}%
\mkern-4mu\cdot\mkern-3mu\cdot\mkern-3mu\cdot\mkern-2mu%
.%
\mkern-4mu\cdot\mkern-3mu\cdot\mkern-3mu\cdot\mkern-2mu%
B_{\alpha_{n}}^{\beta_{n}}\right)
\ \ {=}%
\\
& \ \ \ \ \ \ \ \ \ \ \ \ \ \ \ \ \ \ \ \ \ \ \
{=}\ \
\ \delta_{\beta_{1}\cdots\beta_{n}}^{\alpha_{1}\cdots\alpha_{n}}\nabla
^{\beta_{i}}C\sum_{j}\left(  B_{\alpha_{1}}^{\beta_{1}}%
\mkern-4mu\cdot\mkern-3mu\cdot\mkern-3mu\cdot\mkern-2mu%
\nabla_{\alpha_{i}}B_{\alpha_{j}}^{\beta_{j}}%
\mkern-4mu\cdot\mkern-3mu\cdot\mkern-3mu\cdot\mkern-2mu%
B_{\alpha_{n}}^{\beta_{n}}\right)  .\nonumber
\end{align}
The summation on the right contains $n-1$ terms since the factor
$B_{\alpha_{i}}^{\beta i}$ is absent from the product. According to the
Codazzi equation (\ref{Codazzi}), the tensor $\nabla_{\alpha_{i}}B_{\alpha
_{j}}^{\beta_{j}}$ is symmetric in its lower indices, i.e.%
\begin{equation}
\nabla_{\alpha_{i}}B_{\alpha_{j}}^{\beta_{j}}=\nabla_{\alpha_{j}}B_{\alpha
_{i}}^{\beta_{j}}%
\end{equation}
and therefore vanishes when contracted with the fully skew-symmetric
$\delta_{\beta_{1}\cdots\beta_{n}}^{\alpha_{1}\cdots\alpha_{n}}$. Thus, the
remaining term vanishes, i.e.%
\begin{equation}
\delta_{\beta_{1}\cdots\beta_{n}}^{\alpha_{1}\cdots\alpha_{n}}\nabla
^{\beta_{i}}C\nabla_{\alpha_{i}}\left(  B_{\alpha_{1}}^{\beta_{1}}%
\mkern-4mu\cdot\mkern-3mu\cdot\mkern-3mu\cdot\mkern-2mu%
.%
\mkern-4mu\cdot\mkern-3mu\cdot\mkern-3mu\cdot\mkern-2mu%
B_{\alpha_{n}}^{\beta_{n}}\right)  =0,
\end{equation}
and we have therefore completed the demonstration of equation (\ref{dIB = 0}).

\section{Conclusion}

We have used the tools of the Calculus of Moving Surfaces to show that the
integral of the Gauss-Kronecker curvature of a smooth closed surface is
independent of shape. Of course, Hopf's curvatura integra theorem is a
considerably stronger statement than what we have shown here, since it also
identifies the actual value of the integral in terms of the surface's degree.
Our goal, however, was to showcase the reach of the Calculus of Moving
Surfaces, including its usefulness for problems whose formulations do not, on
their face, involve any moving surfaces at all.

\bibliographystyle{abbrv}
\bibliography{Classics,PGrinfeld,Tensors}

\end{document}